\UseRawInputEncoding
\newcommand{\C}{\mathbb C}
\newcommand{\N}{I\!\!N} 
\newcommand{\R}{I\!\!R}
\newcommand{\Z}{\mathbb Z}
\newcommand{\noo}[1]{\overset {\mbox{%
\lower1pt\hbox{${\scriptscriptstyle o}$}}}n^{\mbox{%
\lower2pt\hbox{$\scriptscriptstyle #1$}}}}
\newcommand{\calK}{\mathcal K}
\newcommand{\Graph}{Gr}
\newcommand{\Dcal}{\mathcal D}

\newcommand{\SHcal}{\mathcal{SH}}

\newcommand{\Lcal}{\mathcal L}

\newcommand{\Cgap}{\mathcal C_{\scriptscriptstyle{\textrm{gap}}}}
\newcommand{\Fcalsa}{\mathcal F^{\scriptscriptstyle{\textrm{sa}}}}

\newcommand{\spfl}{\mathrm{sf}}
\newcommand{\Det}{\mathrm{det}}

\newcommand{\Imm}{\mathrm{Im}}
\newcommand{\Ker}{\mathrm{Ker}}
\newcommand{\im}{\mathrm{Im}}
\newcommand{\Dim}{\mathrm{dim}}
\newcommand{\Codim}{\mathrm{codim}}
\newcommand{\ind}{\mathrm{ind}}

\newcommand{\sgn}{\mathrm{sign}}
\newcommand{\dd}{\mathrm d}
\newcommand{\Id}{\mathrm{Id}}

\newcommand{\fs}{{\Phi_S(H)}}

\newcommand{\la}{\lambda}

\newcommand{\FP}{\mathcal{F} }
\newcommand{\ra}{\rightarrow }
\documentclass[10pt,letterpaper]{amsart}
\usepackage{amssymb}
\usepackage{times}

\def\f{\Phi(X,Y)}
\def\f0{\Phi_{0}(X,Y)}

\def\fsi{\Phi_{S}^{i}(H)}
\def\ab{[a,b]}
\def\01{[0,1]}
\def\!{^{-1}}

\newcommand{\ie}{{\it  i.e.\thinspace{,\thinspace}}}
\newcommand{\cL}{\mathcal{L}}
\newcommand{\spf}{\operatorname{sf}}

\newcommand{\ip}[2]{\langle #1, #2\rangle}

\newcommand{\hfor}{\hbox{ \,for\, }}

\newcommand{\hforall}{\hbox{ \,for\  all\, }}

\newcommand{\sk}{\medskip}

\newcommand{\bt}{\begin{theorem}}
\newcommand{\et}{\end{theorem}}
\newcommand{\bc}{\begin{corollary}}
\newcommand{\ec}{\end{corollary}}
\newcommand{\bp}{\begin{proposition}}
\newcommand{\ep}{\end{proposition}}
\newcommand{\bl}{\begin{lemma}}
\newcommand{\el}{\end{lemma}}
\newcommand{\br}{\begin{remark}}
\newcommand{\er}{\end{remark}}
\newcommand{\bd}{\begin{definition}}
\newcommand{\ed}{\end{definition}}
\newcommand{\be}{\begin{equation}}
\newcommand{\ee}{\end{equation}}

\def\cab{[a,\,b]}
\def\ab{[a,\,b]}
\def\01{[0,\,1]}

\def\R{{\mathbf R}}
\def\C{{\mathbf C}}
\def\N{{\mathbf N}}

\def\Z{{\mathbf Z}}

\def\ep{\epsilon}

\newcommand{\id}{\operatorname{Id}}

\newcommand{\glpath}[1]{\cab\colon \to{\mathcal GL}(H)} 
\newcommand{\glspath}[1]{\cab\colon \to{\mathcal GL}_{S}(H)}

\newcommand{\bdoc}{  \begin{document}  }
\newcommand{\edoc}{  \end{document}  }

\newcommand{\set}[2]{\mathchoice{\left
\{#1\vphantom{#2}\;\right |
\left .\;#2\vphantom{#1}\right 
\}}{\{#1\vphantom{#2}\,|\,#2\vphantom{#1} 
\}} {\{#1\vphantom{#2}\,|\,#2\vphantom{#1} \}} 
{\{#1\vphantom{#2}\,|\,#2\vphantom{#1} \}}}

\def\square{\hfill\parfillskip=0pt {\vrule 
height8pt 
width7pt 
depth0pt}\par\parfillskip=0pt
plus 1fil\medbreak}

\def\fs{\Phi(H)}
\def\fs#1{\Phi_{S}(#1)}
\def\fsi{\Phi_{S}^{i}(H)}

\def\qf#1{Q_{F}(#1)}

\newcommand{\gln}{{\mathcal GL}(\R^n)}

\newcommand{\glsymh}{{\mathcal GL}_{sym}(H)}

\newcommand{\glsymn}{{\mathcal GL}_{sym}{\mathbf R}^n)}

\newcommand{\linsymh}{{\mathcal L}_{(sym}H)}

\newcommand{\linsymn}{{\mathcal L}_{sym}(\R^n)}

\def\fs#1{\Phi_{sym}(#1)}
\def\fsi{\Phi_{sym}^{i}(H)}

\def\fh{{ \mathcal \Phi}(H)}
 \def\fhs{{ \mathcal \Phi}_S(H)}
 
 \newcommand{\qh}{\cab\colon \to {\mathcal Q}(H)} 
 \newcommand{\qhs}{\cab\colon \to{\mathcal Q}_S(H)}

 \newcommand{\lh}{{\mathcal L}(H)}
 \newcommand{\lhs}{{\mathcal L}_S(H)}
 
  \newcommand{\glh}{  {\mathcal GL}(H) } 
 \newcommand{\glhs}{  {\mathcal GL}_{S}(H)  }

\def\clap#1{\hbox to 0pt{\hss#1\hss}}
\def\mathrlap{\mathpalette\mathrlapinternal}
\def\mathrlapinternal#1#2{\rlap{$\mathsurround=0pt#1{#2}$}}
\newcommand{\ubarint}%
{\mathrlap{\intop\limits^{\clap{\rule[-1ex]{0.4em}{\fboxrule}}}}\!\int}
\newcommand{\lbarint}%
{\mathrlap{\intop\limits_{\clap{\rule[1ex]{0.4em}{\fboxrule}}}}\!\int}

\newcommand{\uint}{\ubarint _{\,a}^{\,b} f}
\newcommand{\lint}{\lbarint _{\,a}^{\,b} f}

\newcommand{\coll}[3]{\{#3\}_{#1}^{#2}}

\newcommand{\prodover}[2]{ \Pi_{#1}{#2}  }

\newcommand{\collover}[2]{\{#1\}_{#2} }
\newcommand{\collfrom}[3]{\{#3\}_{#1}^{#2}}

\newcommand{\ifrom}[3]{\int_{#1}^{#2}{#3}}

\newcommand{\ifromd}[4]{\int_{#1}^{#2}#3(#4)\,d#4}

\newcommand{\iover}[2]{\int_{#1}#2}

\newcommand{\ioverd}[3]{\int_{#1}#2(#3)\,d#3}

\newcommand{\ufrom}[3]{\bigcup_{#1}^{#2}#3}

\newcommand{\uover}[2]{\bigcup_{#1}#2}

\newcommand{\interfrom}[3]{\bigcap_{#1}^{#2}#3}

\newcommand{\interover}[2]{\bigcap_{#1}#2}

\newcommand{\sfrom}[3]{ \sum\limits_{#1}^#2 #3}

\newcommand{\sover}[2]{ \sum\limits_{#1}#2}

\newcommand{\calA}{{\mathcal A}}
\newcommand{\calB}{{\mathcal B}}
\newcommand{\calC}{{\mathcal C}}
\newcommand{\calD}{{\mathcal D}}
\newcommand{\calE}{{\mathcal E}}
\newcommand{\calF}{{\mathcal F}}
\newcommand{\calG}{{\mathcal G}}
\newcommand{\calH}{{\mathcal H}}
\newcommand{\calJ}{{\mathcal J}}
\newcommand{\calL}{{\mathcal L}}
\newcommand{\calM}{{\mathcal M}}
\newcommand{\calN}{{\mathcal N}}
\newcommand{\calP}{{\mathcal P}}
\newcommand{\calR}{{\mathcal R}}
\newcommand{\calS}{{\mathcal S}}
\newcommand{\calT}{{\mathcal T}}
\newcommand{\calU}{{\mathcal U}}
\newcommand{\calV}{{\mathcal V}}
\newcommand{\calW}{{\mathcal W}}
\newcommand{\calO}{{\mathcal O}}

\newcommand{\caloh}{{\mathcal O}}

\def\fs {{\Phi}_{S}(H)}

\def\fsi{{\Phi}_{S}^{i}(H)}
\def\qf#1{Q_{F}(#1)}
\def \sp#1#2{\left\langle #1,#2 \right\rangle}
\def\ptq#1{#1\colon [a,b]\to Q_{F}(H)}
\def\sf#1{ \spf(#1,[a,\,b])}

\def\l2{^{1/2}}
\def\dd#1{(#1-\delta,#1+\delta)}

\def\ids{(\id  +|S|)^{1/2}}

\newcommand{\ipfstar}[2]{\langle #1, #2\rangle_{H^{1/2}_*}}

\newcommand{\lo}{\lambda_0}

\newcommand{\no}[2]{\Vert #1 \Vert_{#2)}}
\newcommand{\nof}[2]{\Vert #1 \Vert_{#2^{1/2}}}
\newcommand{\nov}[1]{\Vert #1 \Vert_{V}}
\newcommand{\nog}[2]{\Vert #1 \Vert_{#2}}
\newcommand{\hg}[1]{H_{#1}}
\newcommand{\hf}[1]{H_{#1^{1/2}}}
\newcommand{\sroot}[1]{(I+|#1|)^{1/2}}
\newcommand{\ush}{US(H)}
\newcommand{\Lag}{\Lambda({\bf R}^{2n})}
\newcommand{\m}{\mu_{\operatorname{mas }}}
\newcommand{\llam}{\ell(\lambda)}
\newcommand{\pathl}{\ell\colon \ab \to \Lambda(\R^{2n}}
\newcommand{\elp}{\ell^\prime}
\newcommand{\lamus}{\ell^*}
\newcommand{\lamls}{\ell_*}
\newcommand{\lamo}{\lambda_{0}}
\def\bvpe{H^1_{\ell^\prime,\ell}}
\def\fbvlam{H^{1/2}_{\lambda}}
\def\fbvstar{H^{1/2}_{*}}
\def\bvoe{H^1_{\ell_0,\ell}}
\def\bv{H^1_{\ell_0,\ell_0}}
\def\bvstar{H^1_{*}}
\def\bvlam{H^1_{\lambda}}
\def\bvllam{H^1_{\ell^\prime(\lambda),\ell(\lambda)}}
\def\bvT{H^1_{\lambda}([-T,T])}
\def\l2T{L^2([-T,T])}

\def\l{\lambda}

\def\adpath{$L\colon [a,\,b]\to \Phi_S(H)$ }

\newcommand{\bT }{\bar{T}}
\newcommand{\bP }{\bar{P}}
\newcommand{\bL }{\bar{L}}
\newcommand{\bS }{\bar{S}}
\newcommand{\bare}{\bar{e}}
\newcommand{\ET }{E_-(T)}

\numberwithin{equation}{section}

\title[Spectral flow and variational bifurcation]{Spectral flow and variational bifurcation}

\author[J. Pejsachowicz]{Jacobo Pejsachowicz}
\address{INDAM at Politecnico di Torino, \hfill\break\indent 24 Corso Duca degli
Abruzzi, \hfill\break\indent10129, Torino, TO,
 Italy}

\email{jacobo@polito.it}
\dedicatory{Dedicated to Victor Zvyagin on his 75 years anniversary} 
\thanks{The author is supported by GNAMPA-INDAM}

\subjclass[2020]{58E07, 58J30, 53D12, 47J15, 47A53.}

\begin{document}


\theoremstyle{plain}\newtheorem{teo}{Theorem}[section]
\theoremstyle{plain}\newtheorem{prop}[teo]{Proposition}
\theoremstyle{plain}\newtheorem{lem}[teo]{Lemma}
\theoremstyle{plain}\newtheorem{cor}[teo]{Corollary}
\theoremstyle{definition}\newtheorem{defin}[teo]{Definition}
\theoremstyle{remark}\newtheorem{rem}[teo]{Remark}
\theoremstyle{definition}\newtheorem{example}[teo]{Example}
\theoremstyle{definition}\newtheorem{defins}[teo]{Definitions}
\theoremstyle{plain}\newtheorem*{convention}{Convention}
\theoremstyle{definition}\newtheorem*{defin0}{Definition}

\begin{abstract}
We show that the principle "nonvanishing of  spectral flow of the linearization along the trivial branch entails bifurcation of nontrivial solutions ", proved in \cite{FPR} for critical points of one parameter families of $C^2$ functionals with Fredholm Hessian,  holds true for variational perturbations of paths of unbounded self-adjoint Fredholm operators with a fixed domain. 
\end{abstract}

\maketitle

\tableofcontents




\begin{section}{Introduction} \label{sec:intro}

\noindent{\em A spectre is haunting Europe...} 
 
\vskip5pt

This article is devoted to the proof of a bifurcation theorem for variational perturbations of one-parameter families of unbounded self-adjoint Fredholm operators, generalizing an analogous  bifurcation theorem for one-parameter families of $C^2$ functionals with Fredholm Hessian obtained in \cite{FPR}.  The result is still of the same type: "nonvanishing of the spectral flow of the linearization entails bifurcation of nontrivial solutions", however the setting and the method in proof are completely different from the ones used in \cite{FPR}. The main change pertains to the construction of  the spectral flow. 
 A construction of this well known invariant of a path \cite{APS,Ph} of bounded self-adjoint Fredholm operators  was introduced in \cite{FPR} by reducing a given path of bounded operators to a path of the form  $J+\calK$ where $J$ is a fixed symmetry ($J^2=\id$) and $\calK$ is a path of compact self-adjoint operators. To this aim in \cite{FPR}, we used the cogredience operation $\calA \mapsto \calM^* \calA \calM $ of a path of isomorphisms  $\calM$ (called cogredient parametrix) on a given path of operators   $\calA.$  Applied to the path of Hessians along the trivial branch of critical points of a one-parameter family of  $C^2$ functionals this reduction led to a simple proof of the main bifurcation theorem by Lyapunov-Schmidt method. 
  
In dealing with unbounded operators, the definition of parametrix necessarily needs to be changed. In this setting, a {\em cogredient parametrix} of a path of unbounded operators is a path of bounded operators of the form  $\calM =(\calA + \calK)^{-1},$  where  $\calK$ is a path of compact self-adjoint operators "inverting" the path $\calA.$ 
 However, here we will not define the spectral flow using the above reduction, as  was done in \cite{FPR}.  Instead, we will use a different approach to the spectral flow, devised by Nicolaescu in \cite{ Nic97}.  
 
 The central theme in \cite{ Nic97}  is a systematic use of the symplectic reduction applied to graphs of unbounded self-adjoint operators, viewed as Lagrangian subspaces of the symplectic Hilbert space $H\times H.$ The use of symplectic reduction in this setting resembles the well-known Lyapunov-Schmidt method. Even more suggestive is the fact that the finite-dimensional Lagrangian bundle obtained from the family of graphs of operators by an appropriate symplectic reduction corresponds to the index bundle of a family of general Fredholm operators, used in \cite{P,Wa}  to detect bifurcation of solutions in problems parametrized by compact spaces.  Therefore, the main reason for our choice of Nicolaescu's approach here is that it suggests the possibility of an extension of our results to variational problems depending on several parameters, improving  the results already obtained in \cite{PoWa}.  
Let us point out that, after the appearance of \cite{FPR}, there were other attempts to extend the results of that article to perturbations of unbounded operators (see for example\cite{Wa,FPS}).
Much as this one, none of them appears to be fully satisfactory. Here the problem is as follows: while the spectral flow can be defined for every path continuous in the gap metric, we neither have a proof of our bifurcation theorem in this more general framework,  nor a counterexample that would disprove the possibility of such an extension. 

Below we summarize the content of the article. In the second part of this introduction we present our main results and some corollaries. Section $2$ is a review of the basic theory of symplectic Hilbert spaces.  Here we modify the presentation of the subject by  Nicolaescu  \cite{Nic97} in a form suitable to our purposes. In section $3,$ after reviewing the classical Maslov index,  we introduce its infinite-dimensional counterpart, the generalized Maslov index of a path in the Fredholm Lagrangian Grassmannian.  The spectral flow of a path of unbounded self-adjoint Fredholm operators arises in this framework as the generalized Maslov index of the path defined by  operator's graphs.  It's properties are listed in Propositions \ref{spfl} and \ref{crossing}. Although there are several approaches to the construction of both invariants (see \cite{RS93,RS95,FPR,BLP,CLM,Fu,Wa,W,Nic95,FPR} among others) we believe that Nicolaescu's approach is particularly suitable for nonlinear analysis and bifurcation theory because of the parallelism mentioned above. 

 Section $4$ is devoted to the proof of theorems \ref{thm:riduzione}, and \ref{thm:teoremadibiforca}.
 First, given two paths with values in a finite-dimensional Lagrangian Grassmannian, we construct a path transverse to the suspensions of both  using the analysis of the stable homotopy of this  Grassmannian space, due to F. Latour  \cite{Lat91}.  The construction of the transversal path is strongly related to the existence of a path of self-adjoint operators of finite rank "inverting" a given path $\calA$ of self-adjoint operators. This is the main observation in Theorem \ref{thm:riduzione}. 
 On its turn, the above theorem is used in order to reduce the problem of  existence of bifurcation points of solutions of  \eqref{eq:bifequation} to bifurcation of critical points of a family of  $C^2$-functionals on $H$  whose linearization at the trivial branch is given by a family of compact perturbations of a fixed symmetry $J,$ which was already studied in \cite{FPR}.
   
    The article is partially based on some unfinished notes written in collaboration with A. Portaluri. The author expresses a deep gratitude for his collaboration. 

\end{section}

\smallskip

\subsection{Variational Bifurcation}

We will follow here the setup used in \cite{S} \cite{S1} for the bifurcation from the essential spectrum of unbounded operators. Everywhere here we will use capital calligraphic letters to denote paths. 

Let $\calA= \{A_\la \colon \mathcal D\subset H\to H;  \la \in [a,b]\} $ be a one-parameter family of closed self-adjoint Fredholm operators with a fixed domain $\Dcal.$ That all operators of the family have the same domain will be essential in what follows because  the graph norms    $\|u\|_T = \sqrt{\|u\|^2 + \|Tu\|^2}$  of all closed operators  $T$  with the same domain are equivalent and hence induce on the domain $\Dcal$ the same topology.

 Denoting with $\sigma(T)$  the spectrum of the operator $T,$ we will assume that $0\in \R \setminus \sigma(A_a)\cup \sigma(A_b)$  and hence both end point operators $A_a$ and  $A_b$ have a bounded inverse.     
 Putting  $A = A_a,$  we  introduce the Hilbert space $H_2$  endowing $\Dcal$ with  the graph inner product of $A.$ 

Our next assumption is that  the map  $\la  \ra A_\la $ is continuous with respect to the operator norm in $\cL(H_2;H).$\footnote{While the convergence in gap metric allows to define the spectral flow  it does not appear to be sufficient in order to obtain  bifurcation results analogous to those obtained here.}

 A one parameter family  $\calA \colon \ab \ra \cL(H_2;H)$ verifying the above two assumptions  will be called an {\em admissible  path} in $\cL(H_2;H).$
  By the very definition of its norm, $H_2$ continuously embeds into H.  Therefore, any perturbation $A+\calB$ of a fixed self-adjoint Fredholm operator $A$ by a path $\calB\colon \ab \ra \calL(H)$  of compact self-adjoint operators having invertible endpoints is an admissible path in our sense.   

In order to describe our assumptions about the nonlinear perturbation,  together with $\calD$ we will have to consider its associated form domain.   In our case the absolute value $|A|$ is positive definite and hence admits a positive definite square root. The {\em form domain} $\calD^f$  is the domain of the operator $\sqrt{|A|}.$ We will denote with $H_1$ the  Hilbert space obtained by endowing $\calD^f$ with the graph inner product  of $\sqrt{|A|}.$  It is well known that $H_2$ embeds into $H_1$ as a dense subspace and that for each $\la \in \ab$ the quadratic form $  1/2 \langle A_\la u, u \rangle; u\in \calD $ extends to a bounded quadratic form $q_\la (u)$ on $H_1,$ \cite{S} \cite{S1}.
  \sk 
  
We  will consider a nonlinear perturbation of $ \calA$ given by a Frechet differentiable map 
\be \label{pert} F\colon [a,b] \times H_1 \ra\, H \,\text{such that}\, 
 F_\la(0)=:F(\la,0)=0\,;\, DF_\la(0)=:D_uF(\la,0)=0.\ee

We will also assume that each $F_\la$  admits  a {\em potential}  in the following sense: 
 There exist a $C^2$ function $\psi \colon [a,b] \times H_1 \to R $ such that the differential  of $\psi_\la $ verifies 
\be\label{grad} d\psi_\la (u) [v]= \langle F_\la(u), v\rangle, \hfor\, \text{any}\, u,v \in H_1. \footnote{\ Notice that the scalar product is that of $H.$} \ee

Variational perturbations $F_\la$ of the above type were widely considered in the literature.  They lead to identify the solutions of $A_\la u + F_\la(u) =0$ with the "critical points" of the ill-defined functional $ 1/2 \ip{A_\la u}{u} + \psi_\la u.$   When $H=L^2,$  assuming differentiability of $F_\la $ is unrealistic because it is verified only by affine maps. On the contrary perturbations verifying the above conditions arise in practice as restrictions to the form domain of differentiable maps defined on subspaces $V$ of $L^2$ possessing a stronger topology, but such that $H_1$ embeds continuously in $V.$\, Examples can be found in \cite{AZ, S, Wat}.
 
Under the above  assumptions,  for every $\la \in [a,b],$ the pair  $(\la,0)$ is a solution of the equation
\begin{equation} \label{eq:bifequation} 
A_\la\,u+F(\lambda ,u)=0 
\end{equation} 

This set of solutions is usually called {\em the trivial branch}.  We are interested in finding  points of the trivial branch  from which nontrivial solutions ($u\neq 0$) of   \eqref{eq:bifequation} emerge. 

\begin{defin}  A point  $(\la_*,0)$( or $\lambda_*$)  is a {\em bifurcation
point } for  solutions of the equation  \eqref{eq:bifequation}
from the trivial branch,
 if there exists a sequence
$\{(\lambda_n, u_n)\}_{n \in \N}$   in  $ [a,b] \times \mathcal D$ such that
 $ A_{\la_n} \, u_n \, + \, F(\lambda_n, u_n) \, = \, 0.$
 with  $u_n \, \neq \, 0,$ $\lambda_n \to \lambda_*$ and $u_n \to 0.$
\end{defin}

It is easy to see that  bifurcation can only occur at points $\la$ where $A_\la$ becomes singular, i.e.,  $\ker A_\la\neq \{0\}.$ This provides a necessary condition for bifurcation in terms of the "linearization" $A_\la$ of \eqref{eq:bifequation} along the trivial branch.  Much as in \cite{FPR}, where the bounded case was studied, we associate to every admissible path of self-adjoint Fredholm operators $\calA$  an integral-valued homotopy invariant  $\sf{\calA},$ called  {\em spectral flow}, and prove the following  sufficient condition for bifurcation: 
\vskip5pt

\begin{teo}\label{thm:teoremadibiforca}  Let $ f(\la,u)= A_\la u + F(\la,u)$ be a variational  perturbation of an admissible path  $\calA$  in  $\cL(H_2;H)$  verifying \eqref{pert} and \eqref{grad}.
\begin{itemize} 
 \item [ i)]  If  $\spfl(\calA, [a,b]) \neq \,0$ then  the equation \eqref{eq:bifequation} possesses at least one point  $(\lambda_*,0) \in (a,b)$ of bifurcation of nontrivial solutions  from the trivial branch.
 
\item [ ii)]  Moreover, if  $\,\Ker A_\lambda =\{0\}$ for all but a finite number of $\la$-s,  then \eqref{eq:bifequation}  possesses  at least  $|\spfl(L,I)|/m$ bifurcation points, where 
\[ m=\text{max}\{ \dim \Ker A_\lambda: \lambda \in [a,b]\}.\]  
\end{itemize}
\end{teo}

A {\em nondegenerate  singular point}  of a $C^1$ family $\calA$ is a point $\la_*$ such that  
$\Ker A_{\la_*}\neq {0}$ and the {\em crossing  form} $Q(\la_*,\calA)$ defined on 
$  \Ker A_{\la_*} \,\text{by}\,  Q(\la_*,\calA)(u)= \ip{\frac{d\calA}{d\la}(\la_*)u}{u}$
 is a non degenerate quadratic form. 
\vskip5pt
From the above theorem and Proposition \ref{crossing} it follows: 

\begin{cor}\label{sgnbif} Every nondegenerate singular point $\la_*$ such that $\sgn\, Q(\la_*,A)\neq 0$ is an isolated bifurcation point from the trivial branch.\end{cor}

\begin{cor}\label{sgnbif2}
If either  $\calA$ is $C^1$ and all singular points of $\calA$ are nondegenerate or if $\calA$ is real analytic then the conclusion  $ii)$ of Theorem \ref{thm:teoremadibiforca} holds. 
\end{cor}
\proof  The first part follows from the above corollary. On the other hand, arguing as in \cite{PW}, the set of singular points of an analytic family can be realized as a set of zeroes of a real analytic function, and the second part follows from this.  
\vskip10pt

It can be shown (see Proposition \ref{spfl})  that for  paths $\calA$ of the form $ A_\la =A+ K_\la,$  where  $A$ is a fixed Fredholm  self-adjoint unbounded operator and $K_\la$ compact self-adjoint, the spectral flow depends only on the endpoints of the path and  coincides with the relative Morse index  defined by 
\be\label{remorse} \mu_{rel}(A_a, A_b) = \dim  E^-(A_a)\cap E^+(A_b)- \dim E^-(A_b)\cap E^+(A_a),
\ee 
where $E^{\pm}$ denotes the positive (respectively negative) eigenspace of a self-adjoint operator. 

Therefore, for families  with  linearization along the trivial branch of  the above form, bifurcation arises whenever the relative Morse index $\mu_{rel}(A_a, A_b)\neq 0 .$

If $A_a, \,A_b$ are essentially positive, then the relative Morse index is the difference of the Morse indices of the endpoints 
\be\label{morse}  \mu_{rel}(A_a, A_b) = \mu(A_a)- \mu(A_b)= \dim E^{-}(A_a)-\dim E^{-}(A_b). \ee
 
 Hence we have  
\begin{cor}   
Assume that $ f(\lambda,u) = Au + K_\la u + F(\la, u)$ where $A$ is a self-adjoint positive  Fredholm operator, $K_\la$ compact self-adjoint, and $F$ as above. Then bifurcation arises whenever  $\mu(A + K_a ) \neq \mu(A+K_b).$ 
\end {cor} 

 When the embedding $i\colon \calD(A)\ra H$ is compact  then the above corollary  applied to $f(\la,u) =Au -\la u + F(\la, u)$ takes the form of the Krasnoselkij principle for potential operators: 
 {\em Every eigenvalue $\la $ of $A$ is a bifurcation point for the solutions of the equation $ Au  + F(\la, u)=\la u$  from the trivial branch.}

As we pointed out in the introduction, the proof of  \ref{thm:teoremadibiforca}  is based on a construction of a "cogredient parametrix" analogous to the one defined for bounded operators in \cite{FPR}.  This parametrix is constructed by means of the following auxiliary theorem of independent interest.

 \begin{teo}\label{thm:riduzione} Let $\calA \colon [a,b] \ni
\lambda \mapsto A_\lambda \in \Fcalsa(H)$ be a one-parameter family of unbounded self-adjoint Fredholm operators continuous in the gap topology. Then there exist a
path $K\colon [a,b] \to L(H)$ of bounded self-adjoint operators whose image is contained on a fixed finite-dimensional subspace $F$ of $H$,   such that $A_{\lambda} +K_{\lambda}$ has a bounded inverse for every $\la\in [a,b].$ \footnote{%
In other words, $0$ does not lies in the spectrum of $A_{\lambda}
+K_{\lambda},$ for any $\lambda \in [a,b]$.\/}
\end{teo}

\begin{section}{Symplectic Hilbert spaces}\label{sec:symHilspa}

\begin{defin}
A {\em symplectic Hilbert space\/}({\em sh-space\/})
is a Hilbert space $S$ endowed with a symplectic form $\omega$ and an isometry $J\colon S\rightarrow S$ verifying $\omega (u,v)= \langle Ju,v\rangle.$ 
\end{defin}
\begin{rem}
Due to the skew-symmetry of the symplectic form $\omega,$ any isometry $J$ verifying the above condition  must be a complex structure, i.e.,   $J^{2}=-\Id.$   Using $J$  one can extend the scalars from $\R$ to $\C$ endowing  $S$ with
an extra structure of a Hermitian vector space. \end{rem} 

\begin{rem}

 The notion of symplectic Hilbert space is a very particular case of the scenario considered by Nicolaescu in \cite{Nic97}, where Hilbert modules over Clifford algebras $C^{p,q}$ are studied using techniques typical of the symplectic linear algebra.  The only results  from \cite{Nic97} that we will need in the formulation and proof our bifurcation theorem correspond to the case of modules over the Clifford algebra $C^{0,1}\simeq \C.$ These results are needed for the construction of the spectral flow using  Nicoalescu's approach.  Since in this particular case the proofs of Nicolaescu's results are simpler, we will shortly sketch a few of them in the following sections.
\end{rem}

\begin{example} 
Every Hilbert space having  a bounded  symplectic  form $\omega$  can be  endowed  with  a compatible structure of a symplectic
Hilbert space by eventually modifying the scalar product of the space.

In fact, by Riesz representation theorem there exists  a topological isomorphism $T$  such that 
$\omega (u,v)= \langle Tu,v\rangle.$
Taking the polar decomposition  $T\,=\,JP$ of  $T$ and defining  a new 
scalar product $(u,v)  \equiv \langle Pu,v\rangle$ we get 
\begin{equation*}
\omega(u,v) \,= \,
\langle PJ\, u, v \rangle_1 \, = \, ( Ju, v ),
\end{equation*}
which is a sh-structure with respect to the new scalar product.

In particular any symplectic space $(S,\omega)$
of finite dimension can be endowed with a structure of sh-space. \end{example}

\begin{example}

Let $H$ be any Hilbert space, consider $S(H)=H\times H$
 with the canonical  symplectic
form $\omega((u_1,v_1),(u_2,v_2))= \langle v_2,u_1\rangle -\langle
v_1,u_2\rangle$ coming from the identification of the space $H$ with its dual
$H^*.$ 

Here the associated complex structure $J$ is given by $J(u,v)\, :
=\, (-v,u)$. We will call $S(H)$ the {\em  canonical sh-space
associated with the Hilbert space  $H$}.
\end{example}

Let $S$ be a sh-space.  Given a subspace $W$ of $S,$ we will
denote by $W^{\sharp}$ the orthogonal of $W$ with respect to the
symplectic form $\omega $. In other words
$$W^{\sharp}\,=\,\{v \in S |\quad
\omega(u,v)=0\  \quad \forall u\in W\}.$$
The inner product orthogonal will
be denoted as usual by $W^{\perp}$. The subspaces
$W^{\perp}$ and $W^{\sharp}$ are closed and verify the usual
orthogonality relations. In terms of the complex structure $J$
they are related by 
$ J W^{\sharp}\,=\,  W^{\perp}.$

Let us recall the standard classification of sub-objects in the category of symplectic vector spaces.

\begin{defin}
Let $W$ be a closed subspace of a sh-space $S$. We say that
\begin{enumerate}
\item $W$ is {\em isotropic\/} if $W\subset W^{\sharp}$;
\item $W$ is {\em coisotropic\/} if $ W^{\sharp}\subset W$;
\item $W$ is {\em symplectic\/} if $ W^{\sharp}\cap W= {0}$.
\item $W$ is {\em Lagrangian\/} if $W = W^\sharp.$
\end{enumerate}
In other words, $W$ is Lagrangian if an only if
 $JW = W^\perp$.
 \end{defin}
\begin{example}

In   $S=S(H)$  the
subspaces $ H_0 \, = \, H\times \{ 0\}$ and  $H_1\, = \, \{
0\}\times H$ are always Lagrangian subspaces; moreover  $H_1\,=\, JH_0$.
On the other hand, given a Lagrangian subspace $L$  of a sh-space S, the sh-space
$ S(L)= L\times L $ is isomorphic  in the sh-category  to  $S$ by an symplectic isomorphism sending $L_0$ to $L$ and $L_1$ to $JL.$
\end{example}
\end{section}

\subsection{Fredholm pairs}\label{sec:FreLagGras}

In this subsection, we briefly recall some useful definitions and
 well-known facts about {\em Fredholm pair of subspaces.\/} We will refer to Kato's book \cite[Section IV, 4]{Kat80} in which this definition is given in  more generality.

\begin{defin}\label{def:Fredholmpair}
Let $H$ be a separable Hilbert space and let $V, W$
two closed subspaces. We say that the pair $(V, W)$ is a {\em
Fredholm pair\/}, if the following conditions hold:
\begin{enumerate}
\item $ V + W $ is closed.
\item $\ind (V,W) = \Dim(V\cap W) -\Codim(V+W)<\infty.$
\end{enumerate}
\end{defin}
  The set of all Fredholm pairs in $H$ will be denoted by  $\FP$ The number $\ind (V,W)$ 
 is called the {\em Fredholm index} of the pair.
The relationship of Fredholm pairs with Fredholm operators is as follows:

If $ T \colon  \mathcal D\subset H \to H$  is an unbounded  Fredholm operator then
  $ (H_0, \Graph\, T) $  is a Fredholm pair in $H\times H$.  Moreover $ \ind\,T = \ind\, (H_0, \Graph\, T).$

On the other hand, if  $H$ is a separable Hilbert space and if 
$V, W$ are  two closed subspaces, then $(V,W)$ is a Fredholm pair
if and only if the orthogonal projector onto $W^\perp$  restricted to $V$ 
is a Fredholm operator. Moreover
$\ind (V,W) \, = \, \ind\Big(P_{W^\perp}\big\vert_V \Big).$
\sk

If $V, W$ are Lagrangian
subspaces of a sh-space $S$ such that their sum is closed,  then we have:
\be \label{2.1}\Codim(V+W)=\Dim(V^{\perp}\cap W^{\perp})
 =\Dim(JV\cap JW)=\Dim(V\cap W).\ee

Hence, a pair of  Lagrangian subspaces is
Fredholm  if and only  $L+M$ is closed and $\Dim(L\cap M)< \infty .$
Moreover, by  \eqref{2.1}, the index $\ind (L,M)$ of a Fredholm pair of Lagrangian subspaces
is always $0.$

\subsection{Topology of the Fredholm Lagrangian Grassmannian}
Given a sh-space $S$, let us denote by $\Lambda= \Lambda(S)$
the set of all Lagrangian subspaces $S$. This
set has a natural topology induced by the  {\em gap metric\/}
$d_g(V,W)=\|P_V-P_W\|_{\mathcal{L}(S)},$
where $P_V$ and $P_W$ are the orthogonal projectors onto $V$ and
$W$ respectively.

In order to simplify notations the space $\Lambda(S(H))$  will be denoted by $\Lambda(H)$ in what follows.

 When $H$ is infinite dimensional, the space  $\Lambda (H)$ is contractible (see \cite{Nic97}). However  the closed subspace 
\begin{equation}\label{eq:fredlagrgras}
\Lambda_0\,=\, \Lambda_0(H) \, = \, \{ V \in \Lambda(H): \  (H_0, V) \in \FP\}
\end{equation}
has a nontrivial topology. The space $\Lambda_0(H) $ will be called  the {\em Fredholm Lagrangian Grassmannian\/.}


\subsubsection{Graphs of self-adjoint operators and Lagrangians}

Let $\mathcal{D} \subset H$ be a dense subspace and let
$T:\mathcal{D}\subset H\rightarrow H$ be an unbounded self-adjoint
operator. Let us recall that  $T$  is self-adjoint if
\begin{equation}\label{eq:autoaggiunti}
\mathcal D  = \mathcal D (T) = \mathcal D (T^*)\quad  \text{and}
\langle x,Ty \rangle
\,=\, \langle T x,y\rangle \ \  \forall x,y \in \mathcal{D}.
\end{equation}
In terms of the symplectic structure on $S(H)$ the geometric meaning of condition  \eqref{eq:autoaggiunti} consists precisely on that the graph $G(T)$ of the operator is a closed Lagrangian subspace of $S(H).$ Clearly $G(T)\in \Lambda_0$ whenever $T$ is Fredholm. 

On the other hand, if a subspace  $L \in \Lambda_0$ is such that  $L \cap H_1 \, = \, 0$, then  L  is the graph of a closed operator $T:\mathcal{D}\subset H\rightarrow H.$  Since $L \cap H_1$  is the orthogonal of  $ W=L^\perp + H_0$ it follows that $ W$ is dense in S(H).  Using  
$\mathcal D+H_1= L+H_1= J(L^\perp+H_0),$ the density of $W$ implies that  the domain $\mathcal D$ of $T$ is dense in $H.$ From our previous discussion in this section it follows that $T$ is a self-adjoint Fredholm operator.  If moreover  $L \, + \, H_1$ is a closed subspace of $S(H),$  then  $W= J (L \, + \, H_1)$ is closed  and $ \mathcal{D}= H.$   But then $T$ is bounded by the closed graph theorem. 
Hence  we proved : 
\begin{prop}\cite[Example 3.11]{Nic97}\label{thm:lagrangianigrafici}
A closed subspace  $L$ of $H\times H$ is the graph of a self-adjoint
Fredholm operator if and only if $L \in \Lambda_0$ and $L \cap H_1 \, = \, 0$.
 Moreover, this operator is bounded if and only if  $L \, + \, H_1$ is closed.
\end{prop}



\subsection{Symplectic reduction}\label{subsec:sympred}

Given a symplectic space $(S,\omega )$ and a closed subspace
$I$ of $S$ it is easy to see that the symplectic form $\omega $
descends to a symplectic form on the quotient space $S/I$  only
if  the restriction of $\omega$ to $I$ is non-degenerate.
A partial remedy to this poverty of quotient objects in the category of
symplectic vector spaces is provided by symplectic reduction. This operation
associates to any isotropic subspace $I$ of $S$ the symplectic
space $I^{\sharp}/I$ with the induced symplectic form.

Let us denote by $\SHcal$ the category of sh-spaces with morphisms given by symplectic transformations which are isometries of the corresponding Hilbert structures.
 We will view $I^{\sharp}/I$ as a sub-object of $S$  in this category. Namely, given a
  (closed) isotropic subspace $I$ of the sh-space S, we
will identify the quotient space $I^{\sharp}/I$ with the
orthogonal complement $S_I$ of $I$ in $I^{\sharp}$ and we will
refer to  $S_I$ as the {\em
symplectic reduction of $S$ modulo $I$  in  $\SHcal.$\/} Thus, by definition:

\begin{equation}\label{3.1}
S_I\,:=\, I^{\sharp }\cap I^{\perp}\,=\, J(I^{\perp })\cap
I^{\perp }
\end{equation}


Notice  that in $\SHcal$ we have $S_I \, = \, S_{JI}.$  This
 follows from \eqref{3.1} and the relation between the symplectic and scalar product orthogonality. 

Let us observe also that the orthogonal
$S_I^{\perp}$ of $S_I$ in $S$ is also a sh-space. Moreover, $S$ splits into a direct sum of $S_I$ and $S_I^{\perp},$ that is orthogonal for both the symplectic and Hilbert structure on $S$.
Since,
\begin{equation}\label{eq:ortogonalisimplettici}
S_I^{\perp}\,=\,[I^{\perp}\cap (JI)^{\perp}]
^{\perp}\,=\, I\oplus JI,
\end{equation}
with $I$ and $JI$ transverse Lagrangian subspaces of
$S_I^{\perp},$ it follows  that $S_I^{\perp}$ is
canonically isomorphic to the sh-space $S(I).$

Therefore   in $\SHcal,$ we have
\begin{equation}\label{3.2}
S\, \equiv \,  S_I\oplus S(I).
\end{equation}

Let us recall  that, in the  terminology of  \cite[Section
4]{Nic97}  two closed subspaces $V,W$  have a {\em clean intersection} if
$V\cap W=\{0\}\quad  \text{and}\  V+W\ \text{is closed}. $

If $ I $  is isotropic subspace of $S,$ a  Lagrangian subspace $L$ is said to be { \em clean mod $I$} if $L$ and $I$ have a clean intersection. Notice that if  $I$ itself is Lagrangian
then $L$ has a clean intersection with $I$ if and only in they are transverse. We will denote with  $\Lambda_I(S)$  the subspace of $\Lambda (S)$ whose elements are clean Lagrangians mod $I.$  An isotropic subspace $I$ is called {\it cofinite} if dim$ S_I<\infty$.

In what follows  we collect in
a single proposition some of the results in \cite[Section
4]{Nic97} which we will use here.

\begin{prop} \label{thm:propNico}
 Let $I$ be an isotropic subspace of $S$ and let $L$
be clean mod $I$ then

i) the  image $L_I$ of $L\cap I^\sharp$  by the orthogonal projection onto $S_I$ is a
Lagrangian subspace of $S_I$ called { \it the symplectic reduction of
$L$ modulo $I$.}

ii) $ L_I=(L+I)\cap S_I .$

iii) If $I$ is cofinite, the map
$\rho^I:\Lambda_I(S)\rightarrow \Lambda (S_I)$ assigning to each
clean  mod $I$ Lagrangian  its symplectic  reduction mod $I$ is continuous.

iv) The map $\rho^I$ is a homotopy equivalence between $\Lambda_I(S)$ and $\Lambda (S_I).$  It possesses a right inverse $\varepsilon^I :\Lambda (S_I)\rightarrow \Lambda_I(S) $ defined by
$\varepsilon^I (L)=L+JI $ such that  $ \varepsilon^I \circ \rho^I $ is homotopic to $\Id.$
\end{prop}

\begin{proof}
We refer to \cite{Nic97} section 4 for the proofs in the general
setting of $C^{p,q}$-Hilbert modules. For reader's convenience we
shall sketch them in the only case we will need here.  Namely, when $S=S(H)$, and the isotropic subspace is  $I_0=I\times \{0\}\subset S(H),$ where $I$ is closed subspace of $H$ of finite codimension.

Denoting by $F=I^{\perp}$,
we have that $S_{I_{0}}=I_0^{\perp} \cap JI_0^{\perp }= F\times F,$
with its natural Lagrangian frame $F_0 , F_1.$ Moreover $S_{I_0}^{\perp
}=I\times I$ and therefore the splitting \ref{3.2} becomes $$S(H)=F\times
F\oplus I\times I.$$
In order to simplify our notations, in what follows we  will identify  $I$ with $I_0,$ when no confusion arises. 
 
If $L\in \Lambda_I(S)$ then, by definition  its symplectic reduction
$L_I=\rho^I(L)$ is given by
\begin{equation}\label{3.3}
  L_I=\{(P_F x,y) / \ (x,y) \in L, \, \text{and} \   y\in F\}
\end{equation}
where as usual $P_F$ is the orthogonal projector of $H$ onto $F$.

The assertion  ii) follows easily from \ref{3.3}. In order to check i) and iii),
 without loss of generality we can restrict our attention to
Lagrangians $L \in \Lambda_I(S)$ such that  $L_I$ intersects
transversally $F_1$ i.e. those whose reduction is a graph (see \cite [Proposition 1.4]{Nic97}).

By \ref{3.3},\, $L^I \cap F_1=\{0\}$ is equivalent to $L \cap I\times F=\{0\}$. On the other
hand,  since $L+I_0$ is closed and  dim $F<\infty $, we have that $L+I\times F=S(H)$ and hence, $L$
intersect transversally the Lagrangian subspace $I\times F= J F\times I$ of $S(H)$.

Denoting with $ \tilde{H}$ the subspace $F\times I$ we can identify  $S(\tilde{H})$ with $S(H)$
through a tautological symplectomorphism sending $ \tilde{H}_0$ in $F\times I$
and  $\tilde{H}_1$ into $I\times F.$ Since $L$ has clean intersection with $\tilde{H}_1$  it
follows from proposition \ref{2.1} ii) that $L$ is the graph of a self-adjoint operator
$A_L:\tilde{H}\rightarrow \tilde{H}$. In the decomposition $
\tilde{H}=F\times I \; $ the operator  $A_L$ has the matrix form
\begin{equation} \label{3.4}
\begin{pmatrix}
  B_{L} & C_{L} \\
  C^{\ast} _{L} & D_{L}
\end{pmatrix}
\end{equation}
In particular $B_L^{\ast}=B_L$, since $A_L$ is self-adjoint.  But
it is easy to see that the graph $Gr(B_L) =L_I $ and therefore  $L_I$ is a
Lagrangian subspace of $S(F).$

The form \ref{3.4} proves the continuity of
the reduction $\rho^I$ as well. Indeed, for bounded operators
convergence in gap is equivalent to norm convergence. Thus if  $L^n \rightarrow L^0$ it follows
that  $A_{L^n} \rightarrow A_{L^0}$ and in particular
$B_{L^n}\rightarrow B_{L^0}$ which proves the convergence $L^n_I$
to $L^0_I$.

In order to prove iv) let us first notice that $\rho^I
\varepsilon ^I = id$. On the other hand, by ii),   the map $f =\varepsilon ^I \rho^I
\colon \Lambda_I\rightarrow \Lambda_I$ verifies $f(L)^{\perp}\cap L=\{0\}$, or equivalently
$f(L)\cap JL=\{0\}$, for each $L\in
\Lambda_I$. Since $L$ is clean mod
$I_0$ it follows that $JL+f(L)$ is closed and therefore  $JL+f(L)=S(H)$. By the
same argument as before, modulo a symplectomorphism,  each $f(L)$ is
the graph of a  self-adjoint operator $T_L \colon   L\rightarrow L$.
Now, $h\colon [0,1]\times \Lambda_I\to \Lambda_I$ defined by  $h(t,L)=\text{ graph }(tT_L)$
is a homotopy between $id_{\Lambda_I}$ and $f$.

\end{proof}



We close this subsection by relating  the
symplectic reductions modulo two isotropic
subspaces  $I_1, I_2$  of $S$ such that  $I_2\subset I_1$.

 Let $I=I_1 \cap
I_2^{\perp}$ be the orthogonal complement of $I_2$ in $I_1$. Let
$S^i$ be the symplectic reduction of $S$ modulo $I_i$.
 Since
$I_1^{\sharp}\subset I_2^{\sharp}$, $ \, I $ is an isotropic
subspace of $S^2=I_2^{\sharp} \cap I_2^{\perp}$ moreover its
orthogonal in $S^2$ is given by $$ I^{\perp _2}= S^2 \cap
(I_1^{\perp}+ I_2)=I_2^{\sharp}\cap [I_2^{\perp} \cap (I_1^{\perp}
+ I_2)]= I_2^{\sharp} \cap I_1^{\perp}.$$

Similarly one shows that
the symplectic orthogonal of $I$ in $S^2$ is given by $$I^{\sharp _2}= S^2\cap
I^{\sharp}=I_1^{\sharp}\cap I_2^{\perp}.$$
Thus the
reduction of $S^2$ mod $I$ is
 $S^2_I=I_1 ^{\perp}\cap I_2 ^{\sharp}\cap I_1 ^{\sharp}\cap
 I_2^{\perp}=I_1^{\perp}\cap I_1^{\sharp}=S^1$  In other words we have shown
\begin{equation}\label{3.5}
S_{I_1}=[S_{I_2 }]_I
\end{equation}

Keeping the previous notation, and assuming furthermore that $I_1,
I_2$ are cofinite, we shall consider now the action of the symplectic reduction on a
Lagrangian subspace. The bootstrap argument used above shows
that if $L\in \Lambda (S)$ is clean mod $I_1$ then it is also
clean mod $I_2$. Thus $\Lambda^{I_1} \subset \Lambda^{I_2}$ and
moreover from ii) it follows that $L_{I_2}=\rho^{I_2}(L)$ is clean
mod $I$ (here it is enough to check that $L _{I_2} \cap I=\{0\}$).
Using \ref{3.5} and the orthogonal decomposition of $I_2
^{\sharp}=S^2 \oplus I_2$ we obtain that if $u\in I_2^{\sharp}$
then $P_{S^2}(u) \in I^{\sharp _2}$ if and only if $u\in I_1
^{\sharp }$. From this follows that $$L_{I_2} \cap I^{\sharp
_2}=\{ v\in S^2 / v=P_{S_2}(u), \; u\in L\cap I_1^{\sharp}\}$$
Since $P_{S^1}u=P_{S^1}P_{S^2}u$ we obtain that the reduction of $L_{I_2}$ mod $I$ coincides
with the reduction of $L$ mod $I_1$. In synthesis we have  proved
\begin{prop}\label{R}
Let $I_2 \subset I_1$ be cofinite isotropic subspaces of $S$ and
let $I=I_1\cap I_2^{\perp}$ then $\rho^{I_1} = \rho^I \rho^{I_2}$.
\end{prop}


\section{Generalized Maslov Index and Spectral Flow}\label{sec:maslov}
The notion of generalized Maslov index represents one of the many possible ways to extend to the infinite-dimensional setting the seminal ideas of V.Arnold, which introduced the Maslov Index for closed paths of Lagrangian subspaces of $\R^n$ in his topological proof of the Morse index theorem \cite{Arn67}, \cite{Dui76}. 

From the results in \cite{Nic97}, it follows that the Fredholm Lagrangian Grassmannian $\Lambda_0$ is a classifying space for the functor  ${KO}^1,$ and therefore the homotopy groups of this space can be spelled out of the identities $\pi _q(\Lambda_0) = {KO}^1(S^q).$  In particular we have that $\pi_1(\Lambda _f)$ is isomorphic to $\Z$.  The generalized Maslov index is an explicit construction of this isomorphism by means of symplectic reduction to finite dimensions. With this at hand, the spectral flow of an admissible path of unbounded operators  is  defined via the obvious correspondence between the operator and its graph.   
 
    Let us point out that our terminology regarding this point diverges from the one of Nicolaescu.  Dealing with families of $C^{p,q}$ equivariant operators parametrized by a general parameter space $X,$ he reserves the name  "generalized Maslov index" to  the $K^{p,q}$-theoretic class of the Lagrangian subbundle $E$ of $\R^{2n}$  obtained  from the family of graphs by a symplectic reduction.  In the case $X=S^1,$ the Maslov index of the bundle $E$ is what we call the generalized Maslov index here. Nicolaescu and others \cite{Nic97,Fu,BLP} use the same name for both the Maslov index on $R^{2n}$ and its various extensions to infinite dimensions. However, we have chosen to keep this name only for the classical one, introduced by Arnold.

\subsection{Maslov index in finite dimensions}

Since our generalized index is defined by reduction to finite dimensions, we briefly recall  Arnold's definition of the Maslov index  in this case \ie
$H=\R^n, \ S(H) \simeq \R^{2n}.$

We  begin by defining  the Maslov index of a closed path in
$ \Lambda^n=\Lambda (\R^{2n})$ and then introduce its relative version.

Given any Lagrangian $L\in\Lambda^n$ and an
orthonormal basis of $L,$ there exists a unique unitary operator of
$\C ^n\simeq \R^{2n}$ sending the canonical basis of $ \R^n _0$
into the given basis of $L$ and $\R ^n _1 := i\R^n _0$ into $JL$.
Therefore $U(n)$ acts transitively on $\Lambda(n),$ and moreover the
isotropy subgroup of $\R ^n_0$ is clearly $O(n)$. Thus
$\Lambda(n)\simeq U(n)/O(n)$ and hence the square of the
determinant $\Det^2:U(n)\rightarrow S^1$ factors through a smooth
map $\bar \Det^2\colon\Lambda(n)\rightarrow S^1$.

Given a closed path $l \colon S^1\to\Lambda^n,$ the
{\em Maslov index\/} $m(l)$ is the
{\em winding number\/} of $\bar \Det^2 l \colon S^1 \to S^1 $. In other
words, denoting with  $d\theta $ is the angular form on $S^1,$ 
  \begin{equation*}
m(l ):=\frac{1}{2\pi i} \int_l\bar \Det^{2 \,*}(d\theta ).
\end{equation*}
\begin{rem}\label{rem:iso}
 It follows from the above formula that the cohomology class $\mu$  of the pullback  $ \bar \Det^{2 \,*}( \frac{1}{2\pi i} d\theta)$ is integral, i.e., $ \mu \in  H^1(\Lambda (n),\Z).$ 
 
  The  class $\mu $ in called {\em Maslov class\/} and the  Maslov index  is the evaluation of the Maslov class on the loop $l.$ As a consequence we have  that $\mu\colon \pi_1(\Lambda
^n) \rightarrow \Z$ is an isomorphism.
\end{rem}

\smallskip

Now let us discuss an extension of the Maslov index to non-closed paths. Given a Lagrangian subspace $L$ of $\R^{2n},$ this index counts with multiplicities the number of  points  $t$ at which the path  $l(\la)$ fails to be transverse to $L.$ Alternatively, the  Maslov index can be  interpreted as an intersection index of the path with the one-codimensional subvariety $\Sigma(L)=\{L'| \dim L' \cap L \geq 1\},$ of $\Lambda^n.$  

\begin{rem} In Arnold's terminology  $\Sigma({L)}$  is the {\em train of $L$}. \end{rem}

Notice that the complement  $\Lambda^n-\Sigma (L)$ coincides with the open set $\Lambda_L$ of $\Lambda^n$  of all Lagrangian subspaces having a clean intersection with $L.$  This set is clearly diffeomorphic  to the set $\mathcal L_{sym}(JL)$ of all symmetric endomorphisms of $JL, $ and in particular,  is contractible. 

   Given a path $l\colon [a,b] \to\Lambda^n$ such that the endpoints of $l$ belong to  $\Lambda _{L},$ the {\em Maslov index $m_{L}(l)$  of $l$  relative to $L$\/}
is defined as follows.

Since  $\Lambda_{L}$ is connected,  we can find a path $
\tilde{l}$ joining the endpoints of $l$ and contained
in $\Lambda _{L}$. Let $ \bar{l}$ be the concatenation of these
two paths.  We define
\begin{equation*}
m_{L}(l):=m(\bar l )=\frac{1}{2\pi i}
\int _{\bar l}\bar\det^{2 \, \ast }(d\theta ).
\end{equation*}
From the contractibility of $\Lambda _{L}$ it  follows
that $m_{L_0}(l)$ is well defined and moreover it is invariant by homotopies
of paths keeping endpoints in $\Lambda _{L}.$
Moreover,  by \ref{rem:iso} it follows 
that $m_{L}(l)=0$ if and only if $l$ is homotopic to a path in
$\Lambda _{L}=\Lambda^n -\Sigma(L).$

In the proof of our bifurcation theorem, we will need some other related notions,  the     H\"ormander 4-tuple  index, and  the  signature of a
triple of Lagrangian subspaces.

The {\em H\"ormander index}  $h(L_0, L_1, M_0, M_1)$  of a $4$-tuple of Lagrangian subspaces is defined whenever $L_i$ is transverse to $M_j, \;\ i, j=0,1$  and measures  the variation of the Maslov index of a path with respect under changes of the Lagrangian  $L.$ 

Given two Lagrangian subspaces $L_1,\,L_0$  and a  path $l$ whose
endpoints are transverse to both, the respective Maslov indices
are related by
\begin{equation}\label{4.4}
m_{L_1}(l) - m_{L_0}(l)\,=\,h(L_0, L_1 ,l(0),l(1))
\end{equation}

By definition, the {\em H\"ormander index }$h(L_0,L_1,M_0, M_1)$  is the Maslov
index of the closed path $\gamma $ obtained by going from $M_0$
to $M_1$ in $\Lambda _{L_0}$ and then from $M_1$ to $M_0$ in
$\Lambda _{L_1}$. Alternatively, this index can be defined as:
\begin{equation}\label{4.5}
   h(L_0,L_1,M_0, M_1)=m_{L_0}(q ),
\end{equation}
where $ q $ is any path in  $\Lambda (L_1)$ joining $M_0$ to
$M_1$. The right-hand side of \eqref{4.5} is clearly independent
from the choice of $q $ and the formula \eqref{4.4} follows directly from this.
 
 Notice that $ h(L_0,L_1,M_0,M_1)=0$
whenever $M_0$ and $ M_1$ belong to the same component of $\Lambda _{L_0}
\cap \Lambda _{L_1}$. It can be shown that the converse is also
true. We will use this result here and we refer to  \cite{Dui76}
for the proof.

The {\em signature of a
triple of Lagrangian subspaces\/} is defined for  
triples $(L_0,M, L_1)$ with $M$ transverse to $L_i \; i=0,1$.
 Here we will consider only the case in which also $L_i \; i=0,1$ are transverse between them. In this case  the signature of the triple, $\sgn(L_0,M, L_1,)$ is by definition the
signature of a nondegenerate quadratic form $Q$ on $L_0$ constructed  as follows:

Let $L_0\oplus L_1$ be a Lagrangian decomposition of $S(n):=(\R^{2n},\omega)$;
then $M$ can be seen as the graph of an isomorphism  $A\colon L_0\to L_1$.
Let $Q$ be the quadratic form on $L_0$ defined by $Q(v)=\omega (Av,v),$  then the signature of the triple $ \sgn\, (L_0, M, L_1 )$ is  by definition the signature  of $ Q.$

 In \cite{Hor71} H\"ormander proved that the signature of the Lagrangian triples is
related to the H\"ormander index by the formula:
\begin{equation}\label{4.7}
h(L_0,L_1,M_0,M_1)=\frac12[\sgn \, (L_0,M_1,L_1) -
\sgn (L_0,M_0,L_1)].
\end{equation}

From \eqref{4.7} and the definition of the H\"ormander index it follows  that
$M_0, M_1 $ belong to the same component of
$\Lambda _{L_0} \cap \Lambda _{L_1}$ if and only if 
\begin{equation} \label{eqn:sign} \sgn\, (L_0, M_0,L_1)\, =\, \sgn \,(L_0,M_1, L_1).\end{equation}

\subsection{The generalized Maslov index }\label{genmas}

As before, let $S=S(H)$ and let $\Lambda_0\,$ be
the subset of $\Lambda\, = \, \Lambda (S)$ whose elements are all
Lagrangian subspaces such that $(H_0, L)$ is a Fredholm pair.

Given a path $\Lcal\colon[a,b]\to \Lambda_0$ such that $L_a=\Lcal(a)$
and $ L_b=\Lcal(b)$ are both transverse to  $H_0$, we will define the
{\em generalized Maslov index of $\Lcal$\/}. For this, let us choose any
 subspace $I$ of $H$ of finite codimension, such that $ \Lcal(\la)$ has
a clean intersection with $I_0 \,=\, I\times\{0\}$ for all $t$ in
$[a,b]$. The existence of such a subspace
can be shown easily. Indeed, for a
given $\la_0 \in [a,b]$ it is enough to consider $I^{\la_0}$
orthogonal in $H$ to $ \Lcal(\la_0) \cap H_0.$ Then
$I^{\la_0}$ will have a clean intersection with $ \Lcal(\la)$
for $\la$ in a small enough  neighborhood $U_{\la_0}$ of $\la_0$. Taking
as before a finite cover $U_{\la_1},\dots,U_{\la_n}$ of $[a,b]$ and
taking $I\,:=\,\bigcap_{i=1}^{n} I^{\la_i}$ we get a finite codimensional
subspace such that $ \Lcal(\la)$ has clean intersection with
$I$ for all $\la$ in $[a,b].$ This proves the assertion.

Once $I$ is chosen, let $F\,=\, I^{\perp},$ $\overline S \,=\, S(F)$ and $l(\la) \in
\Lambda(\bar{S})$ be the symplectic reduction of $ \Lcal(\la)$ modulo $I.$

By Proposition \ref{thm:propNico} the path
$l\colon[a,b] \rightarrow \Lambda ( \bar{S})$ is continuous with
end points transverse to $F$ and hence the Maslov index
$m_{F_0}(l)$ is well-defined.
\begin{prop}
$m_{F_0}(l)$ is independent of the choice of $I$.
\end{prop}
\proof We prove the proposition for closed paths and then the
general case will follow from the definition of the relative  Maslov index in finite dimensions.  To shorten the notations we will write $m_F$ for $m_{F_0}.$ Clearly it is enough to consider  two cofinite  subspaces  $I, I'$ of $H$ such that $I \subset I'.$

Let us denote by $F, F'$ the orthogonal of $I, I'$ in $H$ respectively
and let $l\,=\,\rho^{I_0}(\Lcal)$, $l'\,=\,\rho^{I'_0}(\Lcal).$ Denoting by $m_F(l)$ the Maslov index of a closed path $l$ in $S(F),$ we have to show that $m_F(l)=m{F'}(l').$
 By Proposition \ref{R}, $l\,=\,\rho^{I_0''} \circ l' $ where 
$I'' :=I \cap I'^{\perp},$ and by
property (iv) in Proposition \ref{thm:propNico} $l'$
is homotopic to $\varepsilon^{ I_0''} \circ l)$
Hence it is enough to show that
$$m_{F'}(\varepsilon ^{ I''_0}\circ l)=m_{F}(l).$$

Recall that by definition, we have
$\varepsilon^{I''_0} \circ l (\la):=l(\la) + JI''_0$.  But then, under the splitting of $F'\times F'$ into
$F\times F \oplus I'' \times I'',$ the unitary matrix sending $F'
\times \{0\}$ into $l (\la) + JI''_0$ which arises in the definition of the
Maslov index splits into
$$U(\la)=\begin{pmatrix}
U_1(\la) & 0 \\
0 & U_2
\end{pmatrix}$$
where $U_2$ is constant and $U_1 (\la)$ sends $F$ into $l (\la)$.
Thus
\begin{eqnarray*}
m_{F'}(l')= m_{F'}(\varepsilon^I\circ l) &=& 
W\big( \Det^2 U\big)= \\
W\big(\Det^2 U_2) \,\cdot\, W\big(\Det^2 U_1\big)&=& W\big(
\Det^2 U_1\big)\,=\, m_{F}(l),
\end{eqnarray*}
where we have denoted by $W$ the {\em winding number\/}.
This proves that $m_{F}(l)$ is independent of the choice
of $I$ and also gives a simple proof of the invariance of
 the classical Maslov index under symplectic reductions.\qed
\smallskip

With this at hand, we can define the infinite-dimensional version of the ordinary Maslov  index:
\begin{defin}
Given a path $\Lcal\colon[a,b]\ra \Lambda _0$ with 
endpoints transversal to $H_0$, the {\em generalized Maslov index of
$\Lcal$\/} is  $$m_{H_0} ( \Lcal)\,:=\,m_{F_0}(l)$$
where $l$ is the reduction of $\Lcal$ modulo $I_0$ and $F_0 = F
\times \{0\}$ with $F \, = \, I^\perp$.
\end{defin}

The standard properties of the finite-dimensional Maslov index hold true for the generalized Maslov index.

In particular,  it is additive under  concatenation of paths and
 vanishes for paths everywhere transverse to $H_0.$ 
 
  If a homotopy $\calH\colon \ab\times[0,1] \ra \Lambda_0$  is transverse to $H_0$ at $\{a, b\}$ and $I= F^\perp $ is a cofinite subspace of $H$ such that $I_0$ is  transverse to all $\calH(\la,t),$ at $\{a, b\},$  then the symplectic reduction mod-$I_0$ of $\calH(\la,t)$ is  transverse to $F_0$ at $\{a, b\}.$   Hence the invariance of the generalized Maslov index under homotopies of paths in $\Lambda_0$ keeping endpoints transverse to $H_0$ follows from the homotopy invariance property of the classical Maslov index.  
An analogous proof leads to the invariance of the Maslov index under symplectic transformations:

 If $M\colon S(H) \ra S(H')$ is a symplectic transformation such that $M(H_0)=H'_0,$ then, for every admissible path    
$\calL\colon[a,b]\ra \Lambda _0(H),$    
\begin{equation}\label{symp} m_{H'_0} (M\circ\calL)\, =\, m_{H_0}( \calL).\end{equation}
\sk

\subsection{Spectral flow of a path of unbounded self-adjoint Fredholm operators.}

In what follows we will denote by $\Cgap(H)$ the
space of all unbounded self-adjoint Fredholm operators on a real
separable Hilbert space $H$ with the topology induced from the gap topology in $\Lambda $  via the graph correspondence  $A \mapsto Gr A.$

   An {\em admissible path} in $\Cgap(H)$ is a continuous  map  $\calA \colon [a,b] \ra  \Cgap(H)$   such that  $ \ker \,A_a=0 =\ker A_b$. 
  
 By Proposition \ref{thm:lagrangianigrafici} the
map  $\calG \colon [a,b]\to \Lambda_0$ defined by
$\calG(\la)\,:=\, Gr\, A_\la$ 
is a path in $\Lambda_0$  that is  transverse to  $H_0$ at the end points $a,b$.

The {\em spectral flow }of an admissible  path $\calA \colon [a,b] \ra  \Cgap(H)$ is defined by \begin{equation} \label{eq:defsflow}
\spfl(\calA, [a,b])\,:=\,m_{H_0}( \calG).
\end{equation}

\begin{prop}\label{spfl} The integer $\spfl(\calA, [a,b])$ has the following properties:
\vskip5pt
\begin{enumerate}

\item[(i)] {\em Normalization.\/}\  If\, $\Ker\, A_\la\,=\,\{0\}$ for all $\la\in I$ then
$\spfl (\calA,I)= 0.$

\vskip5pt
\item[(ii)] {\em Homotopy Invariance.\/}\   Let $\calH\colon
[0,1]\times [a,b]\to \Cgap(H)$  be a family such that for each $s
\in [0,1],$ the path $ \calH(s,\cdot)$ is
admissible. Then
$\spfl (\calH_0,[a,b])\,=\,\spfl(\calH_1,[a,b]).$
\vskip5pt
\item[(iii)]  {\em Concatenation.\/}\ For an admissible path $\calA\colon [a,b]\to
\Cgap(H)$ on the subintervals $[a, c]$ and $[c,b],$
$\spfl (\calA, [a, b])\,=\,\spfl (\calA, [a, c]) \,+\, \spfl (\calA, [c, b]) .$
\vskip5pt

\item[(iv)]{\em  Compact perturbations.\/} If $\calA\colon [a,b] \to \Cgap(H) $ is of the form $A_{\lambda}= A +K_\la $  with $K_\la\in\calL(H)$ compact, then \begin{equation}\label{eq:indicediMorserel}\spfl(\calA, [a,b]) \,=\, \Dim \big (E_-(A_a)\cap E_+(A_b)\big)-\Dim \big(E_- (A_b)\cap E_+( A_a)\big),\end{equation} where $E_\pm (A)$ denote the positive and negative spectral space of the operator $A$.
 \end{enumerate}
\end{prop}
\vskip5pt
\proof Properties (i)-(iii) follow immediately from the corresponding properties of the generalized  Maslov index. In order to prove (iv), one first verifies the formula by direct calculation for particular paths on  $[0,1]$  of the form $J+tK$ with $J^2=\id $  and  $K$ of finite rank.  Then uses Theorem \ref{thm:riduzione} and homotopy invariance to prove the general formula. 

\begin{rem} 
The right-hand side of the equality stated in formula
\eqref{eq:indicediMorserel} is the {\em relative
Morse index\/} $\mu_{rel}(A_a, A_b)$ \cite{FPR} and the above equality not only allows the calculation of spectral flow in many concrete cases but also provides a good heuristic
interpretation of the spectral flow. Namely, $\spfl(A, [a,b])$
counts the number of negative eigenvalues of $A_a$ that move to
the positive ones of $A_b$ minus the number of negative
eigenvalues of $A_a$ that move to the positive side. 
However, in general, $\spfl(A, [a,b])$ will depend on the homotopy class of the path and not only on the endpoints.
\end{rem}

\begin{rem}  When $\dim H$ is finite, then a simple consequence of (iv) is  
\begin{equation*} \text{(v)}\ \  \spfl(\calA, [a,b]) =\mu_{rel}(A_a, A_b) = \mu\, (A_b) - \mu\, (A_a), \end{equation*}\,
where $\mu \,( A)$ is the number of negative eigenvalues of $A,$ \ie the  Morse index   of the quadratic form represented  by $A.$ 
\sk

It was proved in \cite{CFP} that there is a unique integral invariant for admissible paths of bounded self-adjoint Fredholm operators verifying (i), (ii), (v), and additivity under direct sums.
The results of \cite{CFP}  can be used to prove that for bounded operators the spectral flow defined in this article coincides with the one constructed in \cite{FPR}. Indeed, a well-known homotopy deforms paths of the form  $A_\la\oplus B_\la$ into paths $(B_\la \circ A_\la, \id).$ On the other hand by "abstract nonsense" (see \cite[Theorem  6.8, Chap I]{Sp}  the operation induced by the $H-$group product on the sets of homotopy classes of maps with values in the  $H-$group of bounded self-adjoint Fredholm operators coincides with concatenation. Therefore (ii) and (iii) imply  additivity under direct sums, which shows that both constructions of the spectral flow coincide in the case of bounded operators. 
\end{rem}

Finally, let us discuss the  computation of the spectral flow across a singular point $\la_*$ of an admissible $C^1$ path as the signature of a quadratic form defined on the kernel of the corresponding operator. 

Let   $\la_*$  be an isolated point in the singular set 
$ \Sigma(\calA)\, = \{ \la \in [a,b] \colon 
\Ker \, A_\la \neq 0\}.$   Since, for $\la$ sufficiently close to but different from  $\la_*,$  the operator $A_\la$ is injective, the
{\em spectral flow across the singular point $\la_*$\/},\,
$ \spfl(\calA, \la_*)\, :=\, \lim_{\epsilon \to 0}\spfl(\calA,
[\la-\epsilon, \la+\epsilon])$ is well defined integral number.
  
Assuming that $\calA\colon\ab  \ra  \calL(H_1H)$  is $C^1,$ following  \cite{RS95},
for every $\la_* \in \Sigma(\calA)$ we define the {\em crossing form}
$ Q(\calA, \la_*)\colon \Ker A_{\la_*}\rightarrow\R\,$   by
$\,Q(\calA, \la_*)(h)\,=\, \langle \dot A_{\la_*}\,h \,,\,h\rangle,\,$ where dot stands for the derivative.

\begin{prop}\label{crossing} 
If  $Q(\calA, \la_*)$  is nondegenerate, then  $\la_*\in \Sigma(A)$ is isolated
in the singular set of $A$ and moreover,  denoting with $\sgn \, Q(\calA, \la_*)$ the signature of  $ Q(\calA, \la_*) $,
$$\spfl(\calA,\la_*)\,=\, \sgn \, Q(\calA, \la_*).$$
\end{prop}

\proof  We leave the first assertion to the reader. For the second, consider
  the isotropic subspace $ I_0=  \{0\}\times \im A_\la.$ Then the symplectic reduction of $ S(H)_{I_0}$ is  $\ker A_{\la_*}\times  \ker A_{\la_*}$ and the symplectic reduction of $Gr A_{\la_*}$ is $ \Ker A_{\la_*}\times 0$ which is  transverse  to $0 \times   \Ker A_{\la_*}.$  Hence, the symplectic reduction of the $Gr A_\la$  is the graph of a symmetric transformation $\bar A_\la \colon \ker A_{\la_*} \to \ker A_{\la_*} $ if  $| \la -{\la_*}| < \delta$  for small  enough $\delta$. The Maslov index of the path $\la \to Gr\bar A_\la $  on $[\la_*-\delta, \la_*+\delta]$  is the spectral flow of  $ \bar A$ on this interval. But in finite dimension, the spectral flow is given by the difference between the Morse indices of the endpoints. Under our assumption of nondegenracy of the form $Q(A, \la_*),$   the difference in Morse indices at the endpoints coincides with the signature of $Q(A, \la_*)$ \cite{FPR}.
\section{Proofs}
\subsection {Proof of  Theorem \ref{thm:riduzione}}

To prove this theorem we have to resort to the analysis of the homotopy properties of the finite-dimensional Lagrangian Grassmannian due to F. Latour \cite{Lat91}. His results will be used to establish the existence of a path transverse to the suspension of two paths in $\Lambda(S)$ given a priori. With this at hand, we will complete the proof of Theorem\ref{thm:riduzione} by relating the above result to the existence of a path of self-adjoint operators of finite rank that "inverts" a given path $\calA$ of self-adjoint Fredholm operators.

  By a Theorem of Darboux, we can assume, without loss of generality, that our finite-dimensional symplectic space $S$ is the space $ S(n)= R^{2n},$ is endowed with its  standard symplectic structure. By definition, the $k$-th suspension of $S$ is the space $S\oplus S(2k)$, where $S(2k)= R^{4k},$ as in section $2$. 
 
                       We will work with triples of subspaces $(L_0,M,L_1)$ belonging $\Lambda (S)$ such that $M$ transverse to  $L_i$ for $ i=0,1.$
We will denote  with $\widetilde{\Delta}^{2k}$ the co-diagonal in $S(2k)=\R^{2k} \times \R ^{2k}$ defined  by 
$\widetilde{\Delta}^{2k}=\{(x,y,-x,y)\},$                     
                       
\begin{defin}\label{def:sospensione}
Given a triple of Lagrangian subspaces $(L_0,M,L_1)$  as above, the {\em $k$-th suspension} of the triple $(L_0,M,L_1)$  is the triple of Lagrangian subspaces of $S \oplus S(2k),$ given by $(L_0\oplus \R_0 ^{2k},M \oplus \widetilde{\Delta}^{2k}, L_1 \oplus
\R_1 ^{2k})$  where $\R_0^{2k}:=\R^{2k} \times \{0\},$
$\R_1^{2k}:=\{0\}\times \R^{2k}.$ 
\end{defin}
The choice of the co-diagonal in the definition of suspension is due to the fact that the signature $\sgn(\R_0^{2k}, \tilde{\Delta}^{2k}, \R _1^{2k})=0.$ As a consequence of this choice and the additivity,  the signature of a triple  remains  invariant under suspension. 

\begin{prop}\label{thm:trasversale} 
Given a pair of Lagrangian
 subspaces $L_0, L_1  \in  \R
^{2n}$, there exist a subspace  $M \in\Lambda^n$ transverse to both.
\end{prop}

\proof This is certainly true if $L_0$ is transverse to 
$L_,$ since then $\Lambda _{L_0} \cap \Lambda _{L_1}$ is
diffeomorphic to the space of nondegenerate quadratic forms on
$L_0$. The general case can be recast to this  by using a 
symplectic reduction $S$ of $\R^{2n}$ mod $JI$ where $ I=L_0\cap
L_1.$ In this case for $i= 0,1$, $({L_i})_{JI}\simeq L_i/I$ are two
transversal Lagrangian subspaces of $S$ and hence there exists  a
Lagrangian subspace $\widetilde {M}$ of $S$ transversal to both of
them. Now it is enough to take $M\,=\,\widetilde M + JI.$\qed
\sk \sk

\begin{teo}\label{thm:trans}
  Given  a pair of  Lagrangian
paths $l_0,l_1\colon I\rightarrow
\Lambda(S),$ transversal at the end points: 
\begin{itemize}
\item[i)] There exist $k>0$ and a path $p:I\rightarrow
\Lambda(S\oplus \R^{2k})$ such that for any $t\in I,$ $p(t)$ is
transversal to $l_0(t)\oplus \R _0^k$ and $l_1(t) \oplus \R_1^k$.

\item[ii)] Moreover, for $k$ big enough, any pair of paths $l'_0,l'_1$ in $\Lambda^n$   homotopic to   $l_0$ and  $l_1$ through homotopies with transversal end points verifies the same property. 

\end{itemize}
\end{teo}
\proof  We will prove $i)$ assuming  that  $l_0(t)\equiv \R_0^n$ constant. The general case can be easily recasted to this one. 

By Proposition \ref{thm:trasversale}, for each $t\in [0,1]$ one can find $M(t)
\in\Lambda^n$ such that $M(t)$ intersects transversally $\R_0^n$
and $l_1(t).$

The key point here is the use of suspension in order to define a family such that 
$M(t)$ varies continuously with $t$.

Since, being transverse  is  an open condition  and $[0,1]$ is compact, we can find
a partition $0=t_0<t_1< \dots <t_{N+1}=1$ and Lagrangians
$M_1,\dots,M_{N-1}$ in $\Lambda^n$ such that for all $i \, =
1, \dots , N-1 $, each $M_i$ is
 transversal to $\R _0^{n}$ and to $l_1(t)$
for all $t\in [t_{i}, t_{i+1}]$. (We recall that our Lagrangian
the path has endpoints transversal to $\R_0^n$).

Assuming that we have constructed a path
$p\colon [0,t_m] \to \Lambda^{n+s}$ such that $p(t)$ is transverse
to both $\R^{n+s}_0$ and to $l_1(t)\oplus\R^{s}_1$,
we will show how to extend the path $p$ d to a 
path $\widetilde p$ defined on the interval $[0,t_{m+1}]$
which verifies our transversality requests on the bigger interval. 

Let $ \widetilde M_m \,:=\, M_m \oplus \Delta _s$ where $\Delta
_s$ is the diagonal in $\R ^{2s}.$

Clearly, we have  that $\widetilde M_m$ is transverse 
$\R _0 ^{n+s}$ and $l_1(t) \oplus \R _1^s $ for $t_{m}\leq t
\leq t_{m+1}.$

Now let us denote with $\bar n\,=\,n+s,$ with  $\bar
l(t)\,=\,l_0(t) \oplus\R _1^s$  and consider the integer $d\,=\,
\sgn(\R _0^{\bar n},p(t_m),\bar{l}(t_m)) -\sgn(\R _0^{\bar n
},\widetilde M_m, \bar l(t_m))$.

It is easy to see that $-2\bar n \leq d \leq 2\bar n .$
Taking $k\,=\,2\bar n $ and considering a
quadratic form $Q$ on $\R ^k$ with signature $d$, it is clear
that, if we consider the graph of the associated symmetric operator  $N_d\subset
\R ^{2k}$, then
$\sgn(\R _0^k ,N_d, \R _1^k)\,=\,d.$
From the additivity of the signature it follows that :
\begin{eqnarray*}
& &\sgn\,(\R _0^{\bar n+k}, \widetilde  M_m \oplus N_d , \bar
l(t_m) \oplus \R _1^k) = \sgn(\R _0^{\bar n}
,\widetilde M_m,\bar l(t_m)) +d =\\
& &\sgn\, (\R _0^{\bar{n}}, p(t_m),\bar{l}(t_m))
= \sgn\,(\R _0^{\bar{n}+k},p(t_m ) \oplus \tilde{\Delta}
^k,\bar{l}(t_m) \oplus\R _1^k).
\end{eqnarray*}
From this and \eqref{eqn:sign} it follows that  $\widetilde M_m \oplus N_d$
and $p(t_m)\oplus \widetilde \Delta^k$ belong to the same
connected component  $C$ of $\Lambda_{\R_0 ^{\bar n+k}} \cap
\Lambda_{\bar l(t_m)\oplus \R_1^k}$.

Let $n:I\rightarrow C$ be a path such that $n(0)\,=\, p(t_m)
\oplus \widetilde{\Delta }^k$ and $n(1)=\widetilde M_m\oplus N_d$.
By  transversality, for each $s\,\in\,[0,1],$ we can
find a $\delta_s >0$ and a neighborhood $U_s$ of $s$ such that
$n(s')$ is transverse to  $\bar l_1(t) \oplus \R ^k_1$ for all $t$
such that $t_m\leq t\leq t_m +\delta _s,$ and for all $ s' \in
U_s$. Taking a finite subcover $\{ U_{s_1},\dots, U_{s_r}\}$ of
$[0,1]$ and taking $\delta \,=\, \min \delta _{s_j}$ we obtain
that $n(s)$ is transverse to $ \bar l_1(t)\oplus \R^k_1$ for all
$s\in [0,1]$ and all $t$ in $[t_m, t_m +\delta)$.

Let $\varphi \in C\big([t_m,t_m+\delta ]\big)$ be such that:
$  0\leq \varphi\leq 1 \quad \varphi (t_m)=0,  \varphi(t_m+\delta )=1.$
Define $ \widetilde p(t):[0,t_{m+1}]\rightarrow \Lambda(\bar n
+k)$ by
\begin{equation*}
\widetilde p(t) =\left \{
\begin{array}{lll}
p(t) \oplus \tilde{\Delta }^k & 0\leq t\leq t_m,\\
n(\varphi (t)) & t_m\leq t\leq t_m+\delta ,\\
\widetilde M\oplus N_d & t_m+\delta \leq t \leq t_{m+1}.
\end{array}\right.
\end{equation*}
Now, by construction $ \widetilde p(t)$ is transversal to $\R
_0^{n+s+k}$ and to $l_1(t) \oplus \R _1 ^{n+s+k}$. Repeating word
by word this argument a finite number of times  proves $i).$
\vskip 5pt
 The assertion $ii)$ follows readily from Proposition $I.4.2$ in \cite{Lat91}, which asserts that the direct limit of  transverse triples under suspension is a Serre fibration over  $\Lambda^n\times \Lambda^n$ and therefore possesses the homotopy lifting property. 

 \qed
 \vskip10pt 

Now, let us relate the existence of a path of self-adjoint operators of finite rank inverting a given path $\calA$ of self-adjoint Fredholm operators, with the existence of the path of Lagrangian subspaces transverse to the suspension of a given pair of paths, and show that   Theorem \ref {thm:riduzione} is a consequence of Theorem \ref{thm:trans}.  
 
 To see the relation, let us  look  first  at the case of a single Fredholm self-adjoint operator  $A$. Let $I \subset H$ is a subspace of $H$ of finite
codimension such that $I\cap \Ker \,A\,=\, \{0\}$ and let
$F\,=\,I^{\perp}.$ In order to find a self-adjoint bounded operator $K$ of finite rank  such that $A+K$ has a bounded inverse we  look for operators $K$ of the form $K= I_F\, R\,P_F,$  where $R\in \mathcal L_{sym}(F),$  $P_F$ is the orthogonal projector onto $F,$ and $I_F \colon F \to H$ is the inclusion.
Let us consider  the symplectic reduction of
$\Graph\,A$ modulo $I_0.$ By definition this is the Lagrangian subspace $L$ of $S(F)$ given by
$$L\,=\,\{ (P_F (u), Au)\colon \ \  u \in A^{-1}(F)\}\subset F\times F.$$
Let $M \in \Lambda (F)$ be any Lagrangian subspace transverse to
$F_1:=\{0\} \times F$ and to $L$. Being transverse to $F_1,$
 $M$ is the graph of a symmetric operator $-R \in \mathcal L_{sym}(F)$. 
 We claim that  $K= I_F\, R\,P_F$ inverts $A.$
 Indeed by the closed graph theorem is enough to show that $A+K$ is injective (and hence surjective). 
 But if  $A\,u\,+\,I_F\,R\,P_Fu\,=\,0$ then  $A\,u= -I_FRP_Fu.$ Thus $u \in A^{-1}(F)$ and therefore 
 $(P_F (u), Au)\in L.$  On the other hand $(P_Fu,Au)=(P_Fu,-R\,P_Fu)\in M.$  Being $L$ and $M$ transverse both  $Au$ and $P_Fu$ vanish.  Therefore $u \in \Ker A \cap I$  which, by the choice of $I,$ implies that $u=0.$ This proves the claim for a single operator.

 However, when we consider a path, once the symplectic reduction is done, the
existence of a common transversal holds only up to suspension and
hence we have to take care of this too. 

By the same procedure as at the beginning of Section \ref{genmas} we can find a subspace  $ \widetilde I \subset H$ of finite
codimension such that $\widetilde I  \cap  \Ker A_\lambda
\,=\,\{0\},$ for every $\lambda \in [a,b].$ Taking $\widetilde
F \,:=\,\widetilde I^ \perp$, we have by the orthogonality
relations that $\Imm \, A_\lambda \, +\, \widetilde F\,=\,H$ and
hence the family of subspaces $ \widetilde E_\lambda\, =\,A_\lambda^{-1}(\widetilde F)$ form a finite dimensional (trivial) vector bundle $\widetilde E$ over $[a,b].$

 
 Let   $\widetilde e_\la$ be defined  by $  \widetilde e_\lambda u\,:=\, (P_{\widetilde F }u, A_\lambda u).$ 
By  definition, the symplectic reduction of the family $\Lcal(\lambda)$ of graphs  of $A_ \lambda $ modulo ${\widetilde I}_0$ is the image
    $\widetilde l =\{\widetilde l(\la)\};\la\in [a,b] \}$ of the embedding
   $\widetilde e \colon \widetilde E \rightarrow [a,b] \times S(\widetilde F)$  of $\widetilde E$ as a Lagrangian subbundle  of  $ [a,b] \times S(\widetilde F).$ 
   
By Theorem \ref{thm:trans},  there
exist a suspension $\Sigma^{2k}(\widetilde F):=S(\widetilde F \oplus \R^{k})$ of $S(\widetilde F)$ and a path $\widetilde m \colon [a,b]\rightarrow
\Lambda (\Sigma^{2k}(\widetilde F))$ such that $\widetilde m(\lambda)$ is transverse both  to
$\widetilde l(\lambda)\oplus \R^k _0$ and to $\widetilde F_1 \oplus \R^k_1$.

Now it is enough to choose a subspace $I$ of $ \widetilde I$ of
codimension $k$ and do the symplectic reduction of  $S(H)$ modulo $I_0.$  The result is $S(F)$ where   $F\,=\, I^{\perp}.$ 
 Let  $V= F\cap \widetilde F^{\perp},$ then $V_0= V\times\{0\}= I_0^{\perp} \cap \widetilde I_0.$ By Proposition \ref{R},  $S(\widetilde F)= S(F)_{V_0} $   and therefore,   $ S(F) =S( \widetilde F)\oplus S(V).$ 
Any isomorphism sending an orthonormal basis $u_1,\dots, u_k $ of $V$  to the canonical basis of $\R^k$  extends to an sh-isomorphism $\psi$ between  $S(F)$ and  $\Sigma^{2k}(\widetilde F)$ which sends $F_1$ into $(\widetilde F_1)\oplus \R^{k}_1.$ 

Let $\bar l = \psi^{-1}( \widetilde l\oplus \R^{k}_1).$   Then  $\rho^{V_0} (\bar l)= \widetilde l$ and   $\bar m=  \psi^{-1}( \widetilde m)$ is a path of of Lagrangian subspaces pointwise transverse to  both $\bar l$ and the constant path $F_1.$ But, if $l(\lambda)$ is  the symplectic reduction of  $\Lcal(\lambda)$  modulo $I_0,$  by Proposition \ref{R} again,  we have  $\widetilde l= \rho^{V_0}(l).$  Since  both $l$ and $\bar l$ have  the same reduction  modulo $V_0,$ and since, by property $iv)$ of Proposition \ref{thm:propNico}, $\rho^{V_0}$ is a homotopy equivalence, it follows that  the paths $l$ and $\bar l$ are homotopic.  By $ii)$ of Theorem \ref{thm:trans}, for $k$ big enough the existence of a common transversal depends only on the homotopy class of the path, which allows to conclude that the path  $l$ and the constant path $F_1$ possesses  a common transversal path $m$ to both.   Representing  each  $m(\lambda)$ as a graph of a symmetric operator $R_\lambda$ and defining $K_\lambda = I_F\, R_\lambda\,P_F$ we get a path $\calK$ of finite rank  operators such that, for each $\la,$\,  $\calA (\la)+\calK(\la)$  has a bounded inverse.   Notice that the path  $\calK$ is continuous with respect to the operator norm in $L(H)$ because it consist of operators whose image is contained in a finite-dimensional subspace,  and in finite dimensions the graph topology coincides with the norm topology.  
 
\qed
 \vskip10pt


\subsection {Proof of Theorem \ref{thm:teoremadibiforca}}
 
We will reduce the problem of the existence of bifurcation points for solutions of  \eqref{eq:bifequation} to bifurcation of critical points of a  family  of $C^2$- functionals on whole of $H,$ having a special form.

It is enough to consider the strongly indefinite case, \ie when the  positive and negative eigenspaces of  $A_{\lambda} $ are infinite-dimensional. The general case follows from the above one, taking a Cartesian product with an infinite-dimensional Hilbert space (see \cite[Section 6]{FPR} for details).
Taking a spectral resolution $\{E^\mu _\la\}$ be of $C_\la,$ let  $P^{\pm}_\la = \pm \int_0^{\pm\infty } dE^\mu_\la$ be the projectors on the subspaces $H^\pm_\la$ corresponding to the positive and negative part of the spectrum of $C_\la.$ Being  $\calA\colon [a,b]\ra \calL(H_2:H)$ a continuous path of $A$-bounded perturbations, $\calA$ verifies the hypothesis of theorem 5.12 of  \cite{Kat80}. Thus   $0$ belongs to the spectral gap of this family  and therefore  $P^-_\la = E_\la(0)$ and $P^+_\la = \Id - E_\la(0)$  form a family of projectors  $P^{\pm}_\la$ is continuous in the norm topology of $\calL(H).$  Consequently,  the subspaces $H^\pm_\la$ form two Hilbert bundles $\tilde H^\pm$  over $[a,b]$ such that $\tilde H^+\oplus \tilde H^-$ is isomorphic to the trivial bundle  $ \tilde H =[a,b]\times H.$

 On the other hand,  $\la \ra |C_\la|$ is a family of positive definite operators and the unitary operator   $J_\la$ arising in the polar decomposition $C_\la = J_\la |C_\la|$   of the operator $C_\la$ is the symmetry $ J_\la = P^+_\la - P^-_\la$ (see \cite[Chapter 6, Lemma 2.38  ]{Kat80}). Thus, both $J_\la$ and $|C_\la|^{-1}$ are continuous  paths in $\calL(H)$ and hence so is the path $\la \ra |C_\la|^{-1/2}.$ 
 
 Since $\im  |C_\la|^{-1}= \calD(|C_\la|)=\calD$  
It follows that $\im  |C_\la|^{-1/2} = \calD^f.$ Moreover, each $|C_\la|^{-1/2}$ is an isometry between  $H$ into  $H_1$  endowed  with  the scalar product $\ip{u}{v}_\la =\ip {|C_\la|^{1/2}u}{|C_\la|^{1/2}v},$ which restricts to  a topological  isomorphism between $H_1$ and $H_2$ (see \cite{S1}).

  Let $q_\la(u)$ be the bounded extension of the quadratic form $1/2 \ip{A_\la u,}{u}$  to $H_1$ and let $b_\la(u,v)$ be the associated bilinear form.
  It follows from from our assumptions that $(\la,u)$ is a nontrivial solution  of \eqref{eq:bifequation} if and only if 
   \be\label{wsol} b_\la(u,v)=\ip{F(u)}{v}\hforall v\in H_1.\ee
  
 Indeed, clearly \eqref{eq:bifequation} implies \eqref{wsol}. To see the converse, let us observe that by the density of $\calD$ in $H_1$ on has that  the identity  $b_\la(u,v)=\ip{u}{A_\la v}$ holds for every $ u\in H_1, v\in H_2.$
  
 If  $ u_*\neq 0$ verifies \eqref{wsol}, taking  $z=-F_\la(u_*),$ we get
  \[  \ip{u_*}{A_\la v} =b_\la(u_*,v)= \ip{z}{v} \hforall v\in \calD.\] Since $b(u_*,-)$ is a fortiori a bounded functional on $H,$ being $A$ self-adjoint, it follows  that $u_*\in\calD$ and $ A_\la u_*= -F_\la(u_*).$ Hence   $u_*$ solves \eqref{eq:bifequation}.

\sk

  Now let  $\phi_\la $ be the functional defined on $H_1$ by 
\be\label{function}  \phi_\la(u) = q_\la (u) + \psi_\la (u), \ee  where $\psi_\la$ is the potential of $F_\la$ as in \eqref{grad}.   Writing $A_\la $ as $A_\la= C_\la - K_\la,$ a  short calculation shows that the substitution  $u= |C_\la|^{-1/2}v$  leads to a  $C^2$ functional defined over all of $H$ and having the form  $\phi'_\la(v) =
 \|J_\la v\|^2 - 1/2 \ip{K'_\la v}{ v} + \psi'_\la (v),$ 
 with $K'_\la$ compact self-adjoint and $\psi'_\la(v)= o(\|v\|^2).$

Taking  unitary trivializations of $R_\la^\pm \colon H^\pm_a \ra \tilde H^\pm_\la$ we obtain a  continuous family  of isometries $U_\la ; \ U_\la (z)=(R^+_\la P^+_a(z),R^-_\la P^-_a(z)),$ such that  $J_a= U^{-1}_\la J_\la U_\la.$   
  Now, putting $J= J_a= P^+_a - P^-_a,$ the family  of isometries  $\calM \colon [a,b] \ra  \calL( H; H_1)$ defined by $\calM(\la)z=M_\la z = |C_\la|^{-1/2}U_\la z$  transforms the family of functionals $\{\phi_\la; \la \in \ab\}$  into the  family  $\{\phi''_\la; \la \in \ab\}$ of  $C^2$ functionals of the form 
 \be \label{function2} \phi''_\la (z)=\phi_\la(M_\la z)
 =\|Jz\|^2+1/2 \langle K''_\la z, z\rangle  + \psi''_\la (z), \ee   with  $K''$ and  $\psi''$ having  the same properties as above. 
 
 But now $\nabla  \phi''_\la( z) = Jz +K''_\la z + 0(\|z\|),$   and hence  the path  of linearizations of  $\nabla  \phi''_\la$  along the trivial branch
 $\calB\colon \ab \ra \calL(H_1); \,\calB(\la )=B_\la =J+K''_\la $  is a compact symmetric perturbation of a fixed  symmetry  $J$ which arises in the proof of the bifurcation theorem for critical points of  $C^2$ functionals in \cite{FPR}.
 
By the above discussion, $(\la, z) \mapsto  (\la ,M_\la z)$ sends nontrivial solutions of the equation $\nabla_z  \phi''(\la, z)=0$  to nontrivial solutions of the equation \eqref{eq:bifequation}.   Hence, in order to complete the proof of Theorem \ref{thm:teoremadibiforca} repeating word by word the proof of Theorem 1 \cite[page87]{FPR} with our  new definition of the spectral flow,  we have  to show that $ \spf ( \calB, [a,b]) =  \spf ( \calA, [a,b]).$

In order to prove this assertion,  let us notice that the map $N_\la  \colon H_1\times H_1 \ra H\times H$   defined by $N_\la(u,v) =(M_\la u,  M_\la^{-1}v)$  is a symplectic transformation  between $S(H_1)$ and $S(H)$ which send $ H'=H_1\times\{0\}$ into $H\times\{0\}.$ On the other hand, by the very construction of $M_\la,$ it also sends the graph of $B_\la$ into the graph of $A_\la.$ 

 Let $\calN\colon\ab\times[0,1]\ra\Lambda_0(H\times H)$ be the homotopy  defined by $\calN(\la,t)=N_{t\la}(Gr B_\la).$  Then, by the homotopy invariance  and the invariance of the generalized Maslov index under symplectic transformations \eqref{symp} we obtain: 
\[\spf (\calA,[a,b]) =m_H(\calN_1,[a,b])=m_H(\calN_0,[a,b])=
 m_{H'}(Gr\calB,[a,b])= \spf (\calB,[a,b])\]
 \qed


\begin{thebibliography}{99}

\bibitem{AZ} H. Amann,  E. Zehnder,
{\em Nontrivial solutions for a class of nonresonance
problems and applications
to nonlinear differential equations.\/}
Ann.\ Scuola\ Norm.\ Sup.\ Pisa\ Cl.\ Sci.\ (4), {\bf 7} (1980), no. 4, 539-603.


\bibitem{Arn67} V.I. Arnol'd, {\em On a characteristic class entering in quantization conditions}, Func. Anal. Appl.{\bf 1}(1967),
1-13.

 \bibitem{APS}M. F. Atiyah, V. K. Patodi, and I. M. Singer,{\em  Spectral asymmetry and Riemannian geometry. III}, Math. Proc. Cambridge Philos. Soc., {\bf 79}(1976),71-99.
 
\bibitem{BLP}Bernhelm Booss-Bavnbek, Matthias Lesch and John Phillips
{\em Unbounded Fredholm Operators and Spectral Flow} Canad. J. Math. {\bf 57},(2005) 225-250.

 \bibitem{CLM} S. E. Cappell, R. Lee, E. Y. Miller,{\em On the Maslov index}, Comm.\ Pure Appl.\ Math.\  {\bf 47}  (1994),  no.\ 2, 121-186.

 \bibitem{CFP} E. Ciriza, P. M. Fitzpatrick and J. Pejsachowicz,{\em  Uniqueness of spectral flow} Math. Comput. Modelling {\bf 32} (2000), 1495-1501.

\bibitem{Dui76} J.\ J.\ Duistermaat, {\em On the Morse Index in Variational Calculus}, Adv.\ in Math.\ {\bf 21} (1976), 173--195.

\bibitem{FPR} P.M. Fitzpatrick - J. Pejsachowicz - L. Recht. {\em Spectral flow and bifurcation of critical points of strongly-indefinite functional. Part I. General theory}, J. Functional Analysis {\bf 162} (1999), 52-95.

\bibitem{FPS} P.M.Fitzpatrick, J. Pejsachowicz, C.A. Stuart {\em Spectral Flow for Paths of Unbounded Operators and Bifurcation of Critical Points} Unpublished notes (2006).  

\bibitem{Fu} K.Furutani {\em Fredholm Lagrangian Grassmannian and the Maslov index}
Journal of Geometry and Physics {\bf 5} (2004) 269-331.




\bibitem{Hor71} L.H\" ormander, {\em Fourier Integral
 Operators I}, Acta Math.,
vol.\ 127 (1971), 79-183.


\bibitem{Kat80} T. Kato, {\em Perturbation Theory for Linear Operators}, Grundlehren der Mathematischen Wissenschaften, {\bf 132}, Springer-Verlag, New York/Berlin, 1980.

\bibitem{Lat91} F. Latour {\em Transversales lagrangiennes, périodicité
de Bott et formes generatrices pour une immersion lagrangienne
dans un cotangent}, Ann. Sci. École Norm. Sup. (4) {\bf 24} (1991), no. 1, 3-55.


\bibitem{Nic95} Liviu I. Nicolaescu,
{\em The Maslov index, the spectral flow, and decompositions of
manifolds}, Duke Mathematical Journal {\bf 80}, n. 2 (1995),
485-533.



\bibitem{Nic97} L I. Nicolaescu,
{\em Generalized symplectic geometries and the index of families
of elliptic problems}, Memoirs of the American Mathematical
Society {\bf 128}, n. 609 (1997).

\bibitem{P} J.Pejsachowicz, {\em Bifurcation of Fredholm maps I. The index bundle and bifurcation,} Topological Methods in Nonlinear Analysis.  {\bf 38}, 115 - 168.

\bibitem{PW}J.Pejsachowicz, N.Waterstraat {\em Bifurcation of critical points for continuous families of  $C^2$ functionals of Fredholm type} J. Fixed Point Theory Appl. {\bf 17} (2007) 1-32.

\bibitem{PoWa} A.Portaluri, N.Waterstraat {\em A K-theoretical invariant and bifurcation for
homoclinics of Hamiltonian systems}, J. Fixed Point Theory Appl. {\bf19} (2017) 833-851

\bibitem{RS93} J.\ Robbin, D.\ Salamon, {\em The Maslov Index
for Paths}, Topology {\bf32}, No.\ 4 (1993), 827--844.

\bibitem{RS95} J.\ Robbin, D.\ Salamon, {\em The spectral
flow and the Maslov index}, Bull. London Math. Soc. {\bf 27}
(1995), 1-33.

\bibitem{Sp}E.Spanier, {\em Algebraic Topology} McGraw-Hill (1966)

 \bibitem{S} C.A,Stuart, {\em  Bifurcation for Variational Problems When the Linearisation Has No Eigenvalues}Journal of Functional Analysis {\bf  38}, 169-187. 
 
\bibitem{S1} C.A,Stuart, {\em  Spectrum of a self -adjoint operator and Palais Smale conditions}J. London Math. Soc. (2) {\bf 61} (2000) 581-592.

\bibitem{Ph}J. Phillips, {\em Self-adjoint Fredholm operators and spectral flow} Canad. Math. Bull. {\bf 39} (1996), 460-467.

\bibitem{W} C.Wahl {\em A new topology on the space of unbounded self-adjoint operators, K-theory and spectral flow} C*-Algebras and Elliptic Theory II, Trends Math., Birkhauser, Basel, (2008), 297-309.

\bibitem{Wa} N. Waterstraat, The Index Bundle for Gap-Continuous Families, MorseType Index Theorems and Bifurcation, PhD thesis, Georg-August-Universit\"at
Göttingen, (2011).


\bibitem{Wat} N. Waterstraat, Spectral flow and bifurcation for a class of strongly indefinite elliptic systems, Proc. R. Soc. Edinb., Sect. A {\bf 148} (2018) 1097-1113.
\end{thebibliography}
\end{document}